# ON THE DIVISIBILITY OF CLASS NUMBERS OF QUADRATIC FIELDS AND THE SOLVABILITY OF DIOPHANTINE EQUATIONS

Azizul Hoque[1] & Helen K. Saikia [2]

Department of Mathematics, Gauhati University, Guwahati, India-781014

E-mail: [1] ahoque.ms@gmail.com & [2] hsaikia@yahoo.com

## ABSTRACT

In this paper we provide criteria for the insolvability of the Diophantine equation $x^2 + D = y^n$. This result is then used to determine the class number of the quadratic field $\mathbb{Q}(\sqrt{-D})$. We also determine some criteria for the divisibility of the class number of the quadratic field $\mathbb{Q}(\sqrt{-D})$ and this result is then used to discuss the solvability of the Diophantine equation $x^2 + D = y^n$.

2010 AMS Classification: 11R29; 11D61; 11R41.

## 1. INTRODUCTION

The class number problem of quadratic fields is one of the most intriguing unsolved problems in Number Theory and it has been the object of attention for many years of researchers. Ankeny et al. [1], Chakraborty et al. [5], Kishi et al. [14], Nagel [20], Soundararajan [25], Weinberger [26] and Yu [29] studied the class number problem of quadratic fields. It was proved by Nagel [20] that there are infinitely many quadratic number fields each with class number divisible by a given positive integer. Weinberger [26] showed that for all positive integers $n$, there are infinitely many real quadratic fields each with class number divisible by $n$. In [11], we have proved that there exist infinitely many imaginary quadratic fields whose class numbers are divisible by 3. Recently in [10-13], we have found some useful results on the divisibility of class numbers of real and imaginary quadratic fields. The class numbers of quadratic fields can be used in study of Diophantine equations. On the other hand, the class numbers of quadratic fields can be determined by the solvability of Diophantine equations. Thus there has been considerable attention given to the investigation of relationship between the solvability of Diophantine equations and the class numbers of related quadratic fields. The journey to this investigation has been started in 1853 by Lebesgue [15] in which he discussed the solvability of the Diophatine equation $x^2 + D = y^n$. He proved, using an elementary factorization argument, that this equation has no solution for $D = 1$, except $x = 0$ and $y = 1$. Many special cases of this equation have been considered over the several years, but most of the outcomes for general values of $n$ are of honestly recent origin. Fermat showed, a proof is given in [9], that for $D = 2$ and $n = 3$, this equation has only solution, that is, $x = 5$, $y = 3$. In 1943, Ljunggren [16] generalised Fermat's result and proved that this equation has no solution when $D = 2$ and $x \neq 5$. In 1954, this result was alternatively established by Nagell [22]. In 1923, Nagell [21] proved that this equation has no solutions when $D = 3$. This result was duplicated by Brown [4] in 1975 and subsequently by Cohn [7] in 1993. Nagell [21] also proved that this equation has no solution for $D = 5$. In 1955, Nagell [23] showed that for $D = 4$, this equation has only solution, that is, $x = 2$ and $x = 11$. In 1992, Cohn [6] showed that for $D = 19$, this equation has only solution: $x = 18, y = 7, n = 3$ and $x = 22434, y = 55, n = 5$. Finally, Cohn [8] published a historical survey of this equation in 1993. For any positive integer $D$, Wren [27] in 1973 and Blass [2] in 1976, proved the impossibility of the solutions to this equation when $n = 5$. After a couple of year, Blass and Steiner [3] discussed the insolvability of this equation when $n = 7$.

The primary objective of this paper is to investigate the relationship between solvability of the Diophantine equation $x^2 + D = y^n$, $D > 1$ being an integer and the divisibility of the class number of the imaginary quadratic field $\mathbb{Q}(\sqrt{-D})$. This type of relationships can be found in [17], [18], [19], [24] and [28].

## 2. MAIN RESULTS

In this section we discuss and prove our main results. Throughout this section we consider $K = \mathbb{Q}(\sqrt{-D})$ and by $h(K)$ we denote the class number of the quadratic field $K$.

**Theorem 2.1:** Let $D \equiv 1, 2 \pmod{4}$ be a square-free positive integer and $n > 1$ be an odd integer satisfying $na^{n-1} \not\equiv \pm 1 \pmod{D}$, for some integer $a$. If $\gcd(n, h(K)) = 1$, then the Diophantine equation

$$x^2 + D = y^n \qquad (1)$$

has no solutions.

Proof. Let $y$ be an even integer. Then $y^n \equiv 0 \pmod{4}$ and thus Eq. (1) gives $x^2 \equiv 2, 3 \pmod{4}$.

If $D \equiv 1 \pmod{4}$, then Eq. (1) implies $x$ is odd and thus $x^2 \equiv 1 \pmod{4}$. This is a contradiction.

Again if $D \equiv 2 \pmod{4}$, then Eq. (1) implies $x$ is even and thus $x^2 \equiv 0 \pmod{4}$. This is again a contradiction.

Thus $y$ must be an odd integer.

Suppose $p$ be a prime number such that $p | \gcd(x, y)$, then by Eq. (1), $p | D$ and thus $p = D$. Therefore Eq. (1) implies $D^2 | D$. This is a contradiction. Thus $\gcd(x, y) = 1$.

Suppose $(x_0, y_0)$ be an integral solution to the Eq. (1). Then by the above discussion, $y_0$ is odd and $gcd\,(x_0, y_0) = 1$.

Consider the following factorization in the ring of integers $\mathbb{Z}[\sqrt{-D}]$ of the quadratic field $K$.

$$(x_0 + \sqrt{-D})(x_0 - \sqrt{-D}) = y_0^n \tag{2}$$

If $\mathcal{P}$ is a prime ideal in $\mathbb{Z}[\sqrt{-D}]$ such that $\mathcal{P}$ is a common divisor of the ideals $(x_0 + \sqrt{-D})$ and $(x_0 - \sqrt{-D})$, then $\mathcal{P}|(2x_0)$.

Also Eq. (2) gives $\mathcal{P}|(y_0)$. This implies $\mathcal{P} \nmid (2)$ as $y_0$ is odd, and thus $\mathcal{P}|(x_0)$. This contradicts to the fact that $gcd(x_0, y_0) = 1$. Therefore the ideals $(x_0 + \sqrt{-D})$ and $(x_0 - \sqrt{-D})$ are coprime to each other. Thus we can write

$$(x_0 + \sqrt{-D}) = \mathfrak{a}^n$$

$$(u + \sqrt{-D}) = \mathfrak{b}^n$$

for some ideals $\mathfrak{a}$ and $\mathfrak{b}$ in $\mathbb{Z}[\sqrt{-D}]$.

Since $gcd(n, h(K)) = 1$, therefore $\mathfrak{c}^{h(K)}$ is a principal ideal for any ideal $\mathfrak{c}$ in $\mathbb{Z}[\sqrt{-D}]$, and moreover $\mathfrak{a}^n$ and $\mathfrak{b}^n$ are principal ideals, so that the ideals $\mathfrak{a}$ and $\mathfrak{b}$ are principal. Furthermore, since 1 and $-1$ are the only units in $\mathbb{Z}[\sqrt{-D}]$, thus we have

$$(x_0 + \sqrt{-D}) = (a + b\sqrt{-D})^n$$

for some $a, b \in \mathbb{Z}$.

Comparing imaginary part, we see that

$$1 = \binom{n}{1} a^{n-1} b - \binom{n}{3} a^{n-3} b^3 D + \cdots + (-1)^{\frac{n-1}{2}} b^n D^{\frac{n-1}{2}} \tag{3}$$

Thus $b|1$ and hence $b = \pm 1$.

Now Eq. (3) implies $\pm 1 = \binom{n}{1} a^{n-1} - \binom{n}{3} a^{n-3} D + \cdots + (-1)^{\frac{n-1}{2}} D^{\frac{n-1}{2}}$. This implies $na^{n-1} \equiv \pm 1 \ (mod \ D)$. This contradicts to the hypothesis.

This completes the proof.

As a consequence we have the following result.

**Corollary 2.2**: Let $D \equiv 1, 2 \ (mod \ 4)$ be a square-free positive integer and $p$ be an odd prime satisfying $pa^{p-1} \not\equiv \pm 1 (mod \ D)$, for some integer $a$. If $x^2 + D = y^p$ is has integral solution then $p|h(K)$.

We now fix $y$ as a prime, that is consider the Diophantine equation

$$x^2 + D = p^n \tag{4}$$

where $p$ is a prime and $n >$ is an odd integer. Then the following cases arise:

(a) $D \equiv 0 (mod \ 4)$ if one of the following conditions is satisfied:

    (i)    $x$ is odd and $p \equiv 1 \ (mod \ 4)$.

    (ii)    $x$ is even and $p = 2$.

(b) $D \equiv 1 \ (mod \ 4)$ if $x$ is even and $p \equiv 1 \ (mod \ 4)$.

(c) $D \equiv 2 \ (mod \ 4)$ if $x$ is odd and $p \equiv 3 \ (mod \ 4)$.

(d) $D \equiv 3 \ (mod \ 4)$ if one of the following conditions is satisfied:

    (i)    $x$ is even and $p \equiv 3 \ (mod \ 4)$

    (ii)    $x$ is odd and $p = 2$

We are now in a position to state and prove the following result.

**Theorem 2.3**: Let $n > 1$ be an odd integer, $D$ be a square-free positive integer and $p$ be a prime number satisfying the Eq. (4). The ideal class group of the imaginary quadratic field $\mathbb{Q}(\sqrt{-D})$ has an element of order $n$ if one of the following conditions is satisfied:

(I)   $x$ is an integer and $p \equiv 3 \ (mod\ 4)$ satisfying $x^2 < \frac{p^n}{2}$.

(II)  $x$ is an odd integer and $p = 2$ satisfying $x^2 < 2^{n-1}$.

(III) $x$ is an even integer and $p \equiv 1 \ (mod\ 4)$ satisfying $x^2 < \frac{p^n}{p-1}$.

Proof. Let us first consider the conditions given in (I). Then by condition (c) and condition (i) of (d), we see that either $D \equiv 2 (mod\ 4)$ or $D \equiv 3(mod\ 4)$ according as $x$ is odd or even. In either cases, the ring of integers of $\mathbb{Q}(\sqrt{-D})$ is $\mathbb{Z}[\sqrt{-D}]$.

Consider the following factorization in $\mathbb{Z}[\sqrt{-D}]$:

$$(x + \sqrt{-D})(x - \sqrt{-D}) = p^n \qquad (5)$$

Since $(x + \sqrt{-D})$ and $(x - \sqrt{-D})$ are coprime as ideals in $\mathbb{Z}[\sqrt{-D}]$, we have $(x + \sqrt{-D}) = \mathfrak{a}^n$ and $(x - \sqrt{-D}) = \mathfrak{b}^n$ for some ideals $\mathfrak{a}$ and $\mathfrak{b}$ in $\mathbb{Z}[\sqrt{-D}]$ with $\mathfrak{a}\mathfrak{b} = (p)$. Thus the order of $\mathfrak{a}$ in the ideal class group of $\mathbb{Q}(\sqrt{-D})$ is a divisor of $n$.

Let $\mathfrak{a}^m = (u + v\sqrt{-D})$ for some $u, v \in \mathbb{Z}$. Then

$$p^m = u^2 + v^2 D \qquad (6)$$

If $v = 0$, then $p^m = u^2$. This contradicts to the fact that $n$ is odd. Thus $v \neq 0$ and hence Eq. (6) implies $p^m \geq D$.

Again, $x^2 < \frac{p^n}{2}$ implies $D > \frac{p^n}{2}$. Thus $p^m \geq D > \frac{p^n}{2}$. This leads to a contradiction if $m < n$.

Hence $\mathfrak{a}^n = (u + v\sqrt{-D})$ and $\mathfrak{a}^m$ is not principal for any $m < n$. Thus there is an element of order $n$ in the ideal class group of $\mathbb{Q}(\sqrt{-D})$.

Similarly the result holds if we consider the conditions given in (II).

Finally we consider the conditions given in (III). Then by condition (b), we see that $D \equiv 1 \pmod{4}$. Thus the ring of integers of $\mathbb{Q}(\sqrt{-D})$ is $\mathbb{Z}\left[\frac{1+\sqrt{-D}}{2}\right]$. In the ring $\mathbb{Z}\left[\frac{1+\sqrt{-D}}{2}\right]$, we consider the factorization as given in the Eq. (5). However in this case the ideals $\mathfrak{a}$ and $\mathfrak{b}$ must be in $\mathbb{Z}\left[\frac{1+\sqrt{-D}}{2}\right]$.

Let $\mathfrak{a}^m = \left(\frac{u+v\sqrt{-D}}{2}\right)$ for some $u, v \in \mathbb{Z}$. Then

$$4p^m = u^2 + v^2 D \tag{7}$$

If $v = 0$, then $4p^m = u^2$. This contradicts to the fact that $n$ is odd. Thus $v \neq 0$ and hence Eq. (7) implies $4p^m \geq D$.

Again, $x^2 < \frac{p^n}{p-1}$ implies $D > p^n \left(\frac{p-2}{p-1}\right)$. Thus $4p^m \geq D > p^n \left(\frac{p-2}{p-1}\right)$. This leads to a contradiction if $m < n$. Thus $\mathfrak{a}^n = \left(\frac{u+v\sqrt{-D}}{2}\right)$ and $\mathfrak{a}^m$ is not principal for any $m < n$. Hence there is an element of order $n$ in the ideal class group of $\mathbb{Q}(\sqrt{-D})$.

As a consequence we provide the following criteria on the solvability of the Diophantine Eq. (4).

**Corollary 2.4**: $(x_0, y_0)$ is an integral solution of the Eq. (4), where $n > 1$ is an odd integer and $D$ is a positive square-free integer if one of the following conditions is satisfied:

(i) $x_0$ is an integer and $y_0$ is a prime $\equiv 3 \ (mod \ 4)$ satisfying $x_0^2 < \frac{y_0^n}{2}$.

(ii) $x_0$ is an odd integer and $y_0 = 2$ satisfying $x_0^2 < 2^{n-1}$.

(iii) $x_0$ is an even integer and $y_0 \equiv 1 \ (mod \ 4)$ is a prime satisfying $x_0 < \frac{y_0^n}{y_0 - 1}$.

## ACKNOWLEDGEMENTS

Both Authors acknowledge UGC for financial support.

### REFERENCES

1. N. C. Ankeny, E Artin, S. Chowla, The class number of real quadratic fields, Ann. Math. **56** (1952), 479-493.

2. J. Blass, *A note on Diophantine equation* $y^2 + k = x^5$, Math. Comp. **30** (1976), 638–640.

3. J. Blass and R. Steiner, *On the equation* $y^2 + k = x^7$, Utilitas Math. **13** (1978), 293–297.

4. E. Brown, *Diophantine equations of the form* $x^2 + D = y^n$, J. Reine Angew. Math. **274** (1975), 385–389.

5. K. Chakraborty and R. Murty, *On the number of real quadratic fields with class number divisible by 3*, Proc. Amer. Math. Soc., **131** (2002), 41–44.

6. J. H. E. Cohn, *The diophantine equation* $x^2 + 19 = y^n$, Acta Arith. **61**(2) (1992), 193–197.

7. J. H. E. Cohn, *The Diophantine equation* $x^2 + 3 = y^n$, Glasgow Math. J. **35** (1993), 203–206.

8. J. H. E. Cohn, The Diophantine equation $x^2 + C = y^n$, Acta Arith. **65**(4) (1993), 367–381.

9. L. Euler, *Algebra*, Vol. 2, 2nd edi., J. Johnson and Co., London, 1810.

10. A. Hoque and H. K. Saikia, *A note on quadratic fields with class number divisible by 3*, SeMA J., 2015 (in Press). DOI 10.1007/s40324-015-0051-z.

11. A. Hoque and H. K. Saikia, *A family of imaginary quadratic fields whose class numbers are multiples of three*, J. Taibah Univ. Sci., **9** (2015), 399-402.

12. A. Hoque and H. K. Saikia, *On generalized Mersenne Primes and class-numbers of equivalent quadratic fields and cyclotomic fields*, SeMA J., **67** (2015), 71-75.

13. A. Hoque and H. K. Saikia, *On generalized Mersenne Primes*, SeMA J., **66** (2014), 1-7.

14. Y. Kishi, K. Miyake, *Parametrization of the quadratic fields whose class numbers are divisible by three*, J. Number Theory, **80** (2000), 209-217.

15. V. A. Lebesgue, *Sur l'impossibilité en nombres entiers de l'équation $x^m = y^2 + 1$*, Nouvelles Annales des Mathématiques. **9** (1850) 178.

16. W. Ljunggren, *Über einige Arcustangensgleichungen die auf interessante unbestimmte Gleichungen führen*, Ark. Mat. Astr. Fys. **29A** (1943), no. 13.

17. R. A. Mollin, *Diophantine equations and class numbers,* J. Number Theory, **24** (1986), 7-19.

18. R. A. Mollin, *Class numbers of quadratic fields determined by solvability of Diophantine equations,* Math. Comp. **48** (177), (1987), 233-142.

19. R.A. Mollin *Solutions of Diophantine equations and divisibility of class numbers of complex quadratic fields*, Glasgow Math. J. **38** (1996), 195–197.

20. T. Nagel, Über die Klassenzahl imaginär quadratischer zahlkorper, Abh. Math. Sem. Univ. Hamburg, **1** (1922), 140-150.

21. T. Nagell, *Sur l'impossibilité de quelques équations à deux indéterminées*, Norsk. Mat. Forensings Skrifter, **13** (1923), 65–82.

22. T. Nagell, V*erallgemeinerung eines Fermatschen Satzes*, Arch. Math. **5** (1954), 153–159.

23. T. Nagell, *Contributions to the theory of a category of Diophantine equations of the second degree with two unknowns*, Nova Acta Regiae Soc. Sci. Upsaliensis, **16** (2) (1955).

24. A. Pekin, *On Some Solvability Results of Diophantine Equations and the Class Number of Certain Real Quadratic Fields,* Int. J. Contemp. Math. Sciences, **4** (32), (2009), 1605-1609.

25. K. Soundararajan, *Divisibility of class numbers of imaginary quadratic fields*, J. London Math. Soc., **61** (2000), 681–690.

26. P. J. Weinberger, *Real Quadratic fields with Class number divisible by $n$*, J. Number Theory, **5** (1973), 237-241.

27. B. M. E. Wren, $y^2 + D = x^5$, Eureka, **36** (1973), 37–38.

28. H. Yokoi, *On the Diophantine equation and the class number of real subfields of a cyclotomic field,* Nagoya Math. J. 91 (1983), 151-161.

29. G. Yu, *A note on the divisibility of class numbers of real quadratic fields, J. Number Theory*, **97** (2002), 35–44.